\documentclass{amsart}

\usepackage{amsmath}
\usepackage{amscd}
\usepackage{amssymb}
 
\DeclareMathOperator{\vol}{vol}

\newcommand{\be}{\begin{equation}}
\newcommand{\ee}{\end{equation}}
\newcommand{\mbG}{\mbox{\tiny $G$}}
\newcommand{\mbT}{\mbox{\tiny $T$}}
\newcommand{\mbV}{\mbox{\tiny $V$}}

\newcommand{\cal}{\mathcal}

\newcommand{\fg}{{\mathfrak g}}

\newcommand{\ft}{{\mathfrak t}}

\newtheorem{theorem}{Theorem}[section]

\newtheorem{proposition}{Proposition}[section]
\newtheorem{lemma}{Lemma}[section]

\newenvironment{remark}{\medskip 
\noindent {\bf Remark.}}{\mbox{}}

\newenvironment{note}{\medskip 
\noindent {\em Note.}}{\mbox{}}

\newenvironment{conventions}{\medskip 
\noindent {\bf Conventions.}}{\mbox{}}

\newenvironment{definition}{\medskip 
\noindent {\bf Definition.}}{\mbox{}}

\begin{document}

\title
{Degenerate Chern-Weil Theory and Equivariant Cohomology}
\author{Huai-Dong Cao \& Jian Zhou} 
\address{Department of Mathematics\\
Texas A \& M University\\
College Station, TX 77843} 
\email{cao@math.tamu.edu \\
zhou@math.tamu.edu}
\date{}

\begin{abstract}
We develop a  Chern-Weil theory for compact Lie group  action 
whose generic stabilizers are finite in the framework of 
equivariant cohomology. 
This provides a method of changing an equivariant
closed form within its cohomological class to 
a form more suitable to yield localization results.
This work is motivated by our work \cite{Cao-Zho} on reproving
wall crossing formulas in Seiberg-Witten theory,
where the Lie group is the circle.
As applications, we derive  two
localization formulas of Kalkman type for $G = SU(2)$ or $SO(3)$-actions 
on  compact manifolds with boundary. 
One of the formulas is then used to yield a very simple proof of a localization 
formula due to Jeffrey-Kirwan \cite{Jef-Kir2} in the case of $G = SU(2)$
or $SO(3)$.
\end{abstract}
\maketitle


\footnotetext[2]{The authors are supported in part by NSF }

Throughout this paper,  $G$ will be a compact connected Lie 
group, with $\fg$ as its Lie algebra.
Assume that $G$ acts freely on a smooth manifold $P$. 
Then the quotient map $P \rightarrow P/G = M$ gives
$P$ a structure of principal $G$-bundle. 
The celebrated Chern-Weil theory gives us a homomorphism
\begin{equation}
cw: S(\fg^*)^G \rightarrow H^*(M),
\end{equation}
called the Chern-Weil homomorphism. 
Here $S(\fg^*)^G$ is the algebra of polynomials on $\fg$ which 
is invariant under the adjoint representation of $G$ on $\fg$. 
The Chern-Weil construction uses a connection 
$1$-form $\omega \in (\Omega^1(P) \times \fg)^G$ and its 
curvature $2$-form $\Omega = d \omega + \frac{1}{2} [\omega, \omega]$. 
The  equation $d \Omega = [\Omega, \omega]$ can be used to show 
that for any invariant polynomial $F \in S^n(\fg^*)^G$, 
$F(\Omega)$ is the pullback of a closed form on $M$. 
This defines the homomorphism $(1)$. 
Furthermore, for two connections $\omega^0$ and $\omega^1$ with 
curvatures $\Omega^0$ and $\Omega^1$ respectively, 
there is a canonically defined differential form
$T_{(\omega^0, \omega^1)} F$ on $M$, called the transgression form, such that
$$d T_{(\omega^0, \omega^1)} F = F(\Omega^1) - F(\Omega^0).
$$
Therefore, the Chern-Weil homomorphism is independent of the 
choice of $\omega$. 
We call this Chern's formulation.
Cartan \cite{Car1} presented Weil's formulation,
which we shall review in $\S 1$. 
Through Weil's formulation,
Cartan ($\S 5$ in \cite{Car2}) discovered that the Chern-Weil 
homomorphism can be factored as
$$S(\fg^*)^G \stackrel{\phi}{\rightarrow} H^*_G(P) 
\stackrel{(r^G)_*}{\rightarrow} H^*(M),$$
where $H^*_G(P)$ is the  equivariant cohomology of $P$, 
and $\phi$ is the homomorphism which gives $H^*_G(P)$ 
the structure of an $ H^*(BG) \cong S(\fg^*)^G$-module. 
The homomorphism $(r^G)_*$ is induced from a homomorphism on the 
chain level obtained by a similar Chern-Weil 
construction.

In this paper, we shall generalize the above picture to 
the case that the $G$-action on a smooth manifold $W$ is 
only locally free on a dense open set $W^0 \subset W$. 
Using a connection $\omega$ on $W^0$, and a cut-off function $f$, 
we shall construct homomorphisms
\begin{eqnarray*}
cw_f^G: S(\fg^*)^G \rightarrow H^*_G(W),
\end{eqnarray*}
and
\begin{eqnarray*}
(r_f^G)_*: H^*_G(W) \rightarrow H^*_G(W),
\end{eqnarray*}
such that $cw_f^G = (r_f^G)_* \circ \phi$. 
Here $(r_f^G)_*$ is induced from a homomorphism $r_G^f$
at the chain level in Cartan model for equivariant cohomology.
We shall also construct transgression operator to
show that $cw_f^G$ and $(r_f^G)_*$ are independent of 
the choices of connection $\omega$ and the cut-off function $f$. 
An important observation, pointed out to us by Professor Mich\`{e}le
Verge,
is that when one takes $f \equiv 0$, then our calculation shows that
the homomorphism $(r_f^G)_*$ is the identity map.
The main results of this paper are stated in Theorem 2.1-2.6.
We call these results the degenerate Chern-Weil theory. 
We remark that our approach  corresponds to Chern's formulation.
It depends on calculations by brute force.
It is interesting to find a Weil's formulation,
which might make the argument simpler. 

Even though the results of  this paper 
provide an invariant for non-free group actions
(which is interesting in  its own respect),
 the main motivation is to give a method of choosing
a nice representative for an equivariant cohomological class
to obtain localization results.
At the chain level, for  suitable choice of $\omega$ and $f$,
$r_f^G$ gives us a nice way to change an equaivraint closed form 
$\alpha$ within its equivariant cohomological class to $r_f^G(\alpha)$,
with the following property:
in a neighborhood of the sigular set of the group action,
$r_f^G(\alpha) = \alpha$, outside a larger neighborhood,
$r_f^G(\alpha)$ is the pullback of an ordinary differential form
from the quotient.
This provides a simple explanation  for the localization phenomenon
in equivariant cohomology.
When $\deg (\alpha) = \dim (X)$, one often considers integral $\int_X
\alpha$. But we have
$$\int_X \alpha = \int_X r_f^G(\alpha),$$
by dimension reason. However, $r_f^G(\alpha)$ vanishes outside a
neighborhood
of the singular set of the group  (e.g., at where $f = 1$).
So  the only contribution to the integral 
is from near the singular set.
Localization formula could then 
be obtained by shrinking the support of $1-f$.
This is in the same spirit as the proof of 
the localization formula given
in Berline-Getzler-Vergne \cite{Ber-Get-Ver}.
(It might be possible to reprove their formula along this line.)
A similar argument explains why one can expect localization
formula on manifolds with boundary, 
such as Kalkman's formula \cite{Kal1}.
For details, see \S 3. 
It would be interesting to compare our work with the theory of
singular connections of Harvey-Lawson \cite{Har-Law}
which concerns characteristic  classes and singularities of
vector bundle homomorphisms.
For other methods of obtaining localization formulas, see,
e.g., Atiyah-Bott \cite{Ati-Bot} and Witten \cite{Wit}.

In our earlier work Cao-Zhou \cite{Cao-Zho},  a localization formula for
circle action due to Kalkman \cite{Kal2} is used to
obtain wall crossing formulas in Seiberg-Witten theory
due to Li-Liu \cite{Li-Liu} and Okonek-Teleman \cite{Oko-Tel}. 
An important ingredient in \cite{Cao-Zho} is the construction of
degnerate first Chern class for a circle action. 
The results in this paper are nontrivial generalizations from 
circle group to compact Lie groups.
As explained above, 
the application to localization  formula is the main motivation for  
studying degenerate  Chern-Weil theory.

As illustrations of our localization idea,
we prove two nonabelian localization formulas (Theorem \ref{formula1} and 
Theorem \ref{formula2}) of Kalkman type for $G = SU(2)$ and $SO(3)$.  
Theorem \ref{formula1} should be very useful 
in the  study of various wall crossing phenomenon.
In a sequel \cite{Cao-Zho2}, we apply Theorem \ref{formula1}
to study wall crossing phenomenon in symplectic reduction. 
On the other hand,
though moduli spaces in Donaldson theory are in general noncompact
and our results  do not yet readily apply to the study of wall 
crossing phenomenon of Donaldson invariants, 
we believe suitable modifications should yield some results in
this direction. 
Along the same line,
a localization formula of this type for $U(2)$-action should shed some 
lights on the conjectured 
equivalence of Seiberg-Witten theory and Donaldson theory.
We shall leave such issues for future investigations.
As an application of Theorem \ref{formula2}, 
we shall give a very simple proof of the 
nonabelian localization formula
of Jeffrey-Kirwan \cite{Jef-Kir2} in the case of Hamiltonian
$SU(2)$ or $SO(3)$-actions.

The rest of the paper is organized as follows. 
In $\S 1$ we review the equivariant cohomology and fix some notations. 
The degenerate Chern-Weil theory is presented in $\S 2$. 
In $\S 3$ we prove two nonabelian localization formulas of Kalkman type, Theorem \ref{formula1} and
Theorem \ref{formula2}. 
The application of Theorem \ref{formula2} to symplectic reduction 
 is given in $\S 4$.

{\em Acknowledgements}. 
We would like to thank Professors Reese Harvey, 
Blaine Lawson, Claude LeBrun, 
and Mich\`{e}le Vergne for their interest in this work.
We are especially grateful to Professor Vergne
for her insightful comments and valuable suggestions, as well as for
providing us with the much needed reference \cite{Duf-Kum-Ver}.
An observation of hers (see Theorem  2.5) makes our
results potentially more useful.
We greatly appreciate her help in making this revision possible.
We also thank Professor Kefeng Liu for suggesting us to include 
a discussion of Borel model.
The work in this paper is carried out
during the second author's visit at Texas A\&M University.
He thanks the Department of Mathematics and the Geometry,
Analysis and Topology group for hospitality and financial
support. He also finds the lecture notes from Blaine Lawson's
courses on Chern-Weil theory \cite{Law} extremely useful.

\section{Preliminaries on equivariant cohomology}

We will use two differential geometric models, 
the Weil model and the Cartan model, 
for equivariant cohomology. 
For the sake of completeness, we also discuss  Borel model at the end of
this  section.
We refer the reader to Atiyah-Bott \cite{Ati-Bot}, 
Cartan \cite{Car1, Car2}, 
Berline-Getzler-Vergne \cite{Ber-Get-Ver},
Duflo-Kumar-Vergne \cite{Duf-Kum-Ver},
Kalkman \cite{Kal1},
Lawson \cite{Law} and Mathai-Quillen \cite{Mat-Qui} 
and the references therein for more details.

\subsection{Weil algebra}
The Weil algebra \cite{Car1} is the Hopf algebra
$$W(\fg) = \Lambda(\fg^*) \otimes S(\fg^*),$$
where elements in $\Lambda^1(\fg^*)$ have degree $1$, 
and elements in $S^1(\fg^*)$ have degree $2$. 
Let $\{ \xi_a \}$ be a basis of $\fg$, such that
$$[\xi_a, \xi_b] = f_{ab}^c \xi_c, $$
where $f^c_{ab}$'s are the structure constants. 
Let $\{ \theta^a \}$ be a dual basis in $\Lambda^1(\fg^*)$, and 
$\{ \Theta^a \}$ a dual basis in $S^1(\fg^*)$. 
Define the Weil differential $d_w: W(\fg) \rightarrow W(\fg)$ 
by setting
\begin{eqnarray*} 
d_w \theta^a & = &  -\frac{1}{2}f_{bc}^a\theta^b\theta^c +  
	\Theta^a, \\
d_w \Theta^a & = & - f_{bc}^a\theta^b\Theta^c
\end{eqnarray*}
and extending it as a derivation of degree $1$.
There are also contractions $i_a$ and Lie derivatives $L_a$ 
on $W(\fg)$ defined by
$$ \begin{array}{lcl}
i_a \theta^b = \delta_a^b, &  
	& L_a \theta^b = - f_{ac}^b\theta^c, \\
i_a \Theta^b = 0, & & L_a \Theta^b = - f_{ac}^b\Theta^c.
\end{array} $$
Notice that $G$ acts on $W(\fg)$ by extending the co-adjoint 
representation.
Its linearization can be identified with  $L_a$'s. 
It is easy to verify the homotopy formula
$$L_a = d_w i_a + i_a d_w.$$

\subsection{Algebras with Weil structures} We need the following

\begin{definition}  
An {\em algebra with Weil structure} over $G$ is a 
graded differential algebra 
$$({\cal A}^* =
\oplus_{j = 0}^{+\infty}  {\cal A}^j, d)$$
over ${\Bbb R}$, with a left representation 
$$L: G \rightarrow Aut( {\cal A}^*, d)$$
of degree $0$, and a $G$-equivariant linear map $i: \fg
\rightarrow End {\cal A}^*$ of degree $-1$, such that
\begin{enumerate}
\item $i_{\xi} i_{\eta} + i_{\eta}i_{\xi} = 0$, for 
$\xi, \eta \in \fg$;
\item ${\cal L}_{\xi} = d \circ i_{\xi} + i_{\xi} \circ d$, 
where 
${\cal L}: \fg \rightarrow Der({\cal A}^*, d)$
is the linearization of the representation 
$L: G \rightarrow Aut({\cal A}^*, d)$.
\end{enumerate}
\end{definition}

A simple example of algebra with Weil structure is $\Omega(X)$ 
with ordinary contractions, Lie derivatives and the exterior 
differential, for a $G$-manifold $X$. 
Another example is the  Weil algebra. 
Now, given an algebra $({\cal A}^*, d, i, L)$ with Weil 
structure over $G$, define the basic subalgebra by 
$${\cal A}^*_{\mbox{basic}} = 
	\{ \phi \in ({\cal A}^*)^G | i_{\xi} \phi = 0, 
	\forall \xi \in \fg \}.$$
It is straightforward to verify the following

\begin{lemma}
$ {\cal A}^*_{\mbox{basic}}$ is a graded differential
subalgebra of ${\cal A}^*$.
\end{lemma}

The cohomology of the basic subalgebra is called the 
basic cohomology, and denoted by $H^*({\cal A})_{\mbox{basic}}$.
It is well-known that the  basic cohomology  of $W(\fg)$ is $S(\fg^*)^G$.
When $\pi: P \rightarrow M$,  the  basic cohomology of  $\Omega(P)$
is the de Rham cohomology $H^*(M)$.

\begin{lemma}
If $\rho: {\cal A}^* \rightarrow {\cal B}^*$  is a homomorphism 
of algebras with Weil structures over $G$, 
then $\rho$ induces a homomorphism
$$\rho_{\mbox{basic}}: ({\cal A}^*_{\mbox{basic}}, d) \rightarrow 
({\cal B}^*_{\mbox{basic}}, d),$$
and therefore, a homomorphism
$$(\rho_{\mbox{basic}})_*: H^*({\cal A})_{\mbox{basic}} 
\rightarrow H^*({\cal B})_{\mbox{basic}}.$$
\end{lemma}

The above definition and lemmas about algebras with 
Weil structures are taken
from Lawson \cite{Law}. 
They  appeared in Cartan \cite{Car1} and 
Kamber-Tondeur \cite{Kam-Ton} with different terminologies.

If $P \rightarrow M$ is a principal $G$-bundle, a connection 
$\omega = \xi_a \omega^a$ with curvature
$\Omega = \xi_a \Omega^a$
defines a homomorphism of algebras with Weil structures 
\begin{equation}
W(\fg) \rightarrow \Omega(P),
\end{equation}
by sending $\theta^a$ to $\omega^a$ and 
$\Theta^a$ to $\Omega^a$. Applying Lemma 1.2, 
one gets the Chern-Weil homomorphism
$$cw: S(\fg^*)^G \rightarrow H^*(M)$$
by identifying the basic cohomology of $\Omega(P)$ 
with $H^*(M)$. 
We call this way of getting the Chern-Weil homomorphism 
Weil's formulation.
See Cartan \cite{Car1}.
Notice that we can factor $(2)$ as a composition of 
two homomorphisms (see cartan \cite{Car2})
$$W(\fg) \hookrightarrow W(\fg)\otimes \Omega(P) 
\stackrel{cw^p}{\rightarrow} \Omega(P),$$
where the first one is the inclusion, and $cw^p$ is defined 
by extending $(2)$. One can show that $(cw^p_{\mbox{basic}})_*$ is 
an isomorphism.
Indeed, if  
$i^w: \Omega(P) \hookrightarrow W(\fg) \otimes \Omega(P)$ is 
the inclusion, then $cw^p \circ i^w = 1$ on $\Omega(P)$ 
implies that $(i^w_{\mbox{basic}})_*$
is injective. 
Cartan's proof to Theorem 3 in \cite{Car2} shows that
$(i^p_{\mbox{basic}})_*$ is also surjective, and hence an inverse to 
$(cw^p_{\mbox{basic}})_*$.

\subsection{Equivariant cohomology: Weil model and Cartan model}

Let $X$ be a compact smooth $G$-manifold. 
The $G$-action on $X$ induces a homomorphism
from the Lie algebra $\fg$ to the Lie algebra of 
vector fields on $X$. 
Denote by $\iota_a$ and ${\cal L}_a$ the contraction and 
the Lie derivative by the vector field corresponding 
to $\xi_a \in \fg$ respectively.
Consider the tensor product of algebras with Weil structures
$$W(\fg) \otimes \Omega(X),$$
where one uses the diagonal $G$-action, and the contraction 
$i_a \otimes 1 + 1 \otimes \iota_a$,
the Lie derivative  $L_a \otimes 1 + 1 \otimes {\cal L}_a$,
and differential
$d_w \otimes 1 + 1 \otimes d$. 
The corresponding basic cohomology is called
equivariant cohomology (via Weil model), 
and is denoted by $H_G^*(X)$.

Motivated by  the work of Cartan \cite{Car2},
one can also consider the  Cartan model 
which is given by the complex $(\Omega_G(X), D_G)$, 
where $\Omega_G(X) =(S(\fg^*) \otimes \Omega(X))^G$, and
$D_G = 1 \otimes d - \Theta^a \otimes \iota_a$, 
called the Cartan differential. 
When there is only one Lie group involved, 
we will use $D$ for $D_G$. 
Since $D$ is a $G$-invariant operator on 
$S(\fg^*) \otimes \Omega(X)$, 
it then maps $\Omega_G(X)$ to itself. 
Furthermore,
since $\Theta^a \otimes L_a$ acts as zero on $S(\fg^*)$,
we have
$$D^2 = - \Theta^a \otimes {\cal L}_a 
= -\Theta^a  (L_a \otimes 1 + 1 \otimes {\cal L}_a).$$ 
Therefore, $D^2 = 0$ on $\Omega_G(X) = 
(S(\fg^*) \otimes \Omega(X))^G$.

It is possible to identify $H^*_G(X)$ with 
$H^*(\Omega_G(X), D)$ through an isomorphism
$\Psi: W(\fg) \otimes \Omega(X) \rightarrow 
	W(\fg) \otimes \Omega(X) $
defined by
$$\Psi  = \prod_a \exp (- \theta^a \otimes \iota_a) 
	= \prod_a ( 1 - \theta^a \otimes \iota_a).$$
In fact, $\Psi^{-1}$ maps $((W(\fg) \otimes \Omega(X))_{\mbox{\small basic}}, 
d_w \otimes 1 + 1 \otimes d)$ to 
$((S(\fg^*) \otimes \Omega(X))^G, D)$. 
See Cartan \cite{Car2}, Mathai-Quillen \cite{Mat-Qui} 
and Kalkman \cite{Kal1} for more details.

In the case of a principal $G$-bundle $\pi: P \rightarrow M$, 
one can define a homomorphism
$$r^G: S(\fg^*) \otimes \Omega(P) \rightarrow \Omega(P)$$
by $r^G = cw^p \circ \Psi$. It is easy to see that
if $\alpha \in (S(\fg^*) \otimes \Omega(P))^G$ 
and $D \alpha = 0$, then  $d r^G(\alpha) = 0$, 
and $r^G(\alpha)$ is the pullback of a form on $M = P/G$.
As the Cartan model version of the fact that $(cw_{mb}^p)_*$ is an 
isomorphism, $r^G$ induces an isomormphism between
$H^*_G(P)$ with $H^*(M)$.
For a proof, see Duflo-Kumar-Vergne \cite{Duf-Kum-Ver}.
(We thank Mich\`{e}le  Vergne for bringing our attention to
this reference.)

\begin{remark}
When the group action is locally free, i.e., 
all the isotropy subgroups are discrete, 
then $M/G$ is an orbifold \cite{Sat}. 
The above discussions 
carry through if one uses de Rham theory for orbifolds. 
\end{remark}

\subsection{Reduction to the maximal torus}
Let $T$ be a maximal torus of $G$, with Lie algebra $\ft$.
The inclusion $\ft \hookrightarrow \fg$ induces a 
map $\fg^* \rightarrow \ft^*$. Alternatively, if we endow
$\fg$ with a $G$-invariant inner product, then one gets an orthogonal
projection $\fg^* \rightarrow \ft^*$, which can be identified
with the map above. This can be extended to a projection
$p_1: S(\fg^*) \rightarrow S(\ft^*)$. 
Similarly, if we endow $X$ with a $G$-invariant
Riemannian metric, it then induces an inner product on
$\Omega(X)$. So we get a projection $p_2: \Omega(X) 
\rightarrow \Omega(X)^{\mbT}$. Put $p_1$ and $p_2$ together, we get
a projection
$$p_1 \otimes p_2: S(\fg^*) \otimes \Omega(X) 
\rightarrow S(\ft^*) \otimes \Omega(X)^{\mbT},$$
which induces a projection 
\begin{equation}
p: \Omega_{\mbG}(X) \rightarrow \Omega_{\mbT}(X).
\end{equation}
It is an easy exercise to see that $p D_{\mbG} = D_{\mbT} p$, 
hence $p$ induces a homomorphism 
$H^*_{\mbG}(X) \rightarrow H^*_{\mbT}(X)$. 
Let $W$ be the Weyl group,  
then $p$ induces an isomorphism 
$H_{\mbG}^*(X) \cong H^*_{\mbT}(X)^W$. 
For a proof, see e.g. Hsiang \cite{Hsi} 
or Duflo-Kumar-Vergne \cite{Duf-Kum-Ver}.

\begin{note}
In an earlier version, we falsely claim that the image of
$p(\omega_G(X)) = \Omega_T(X)^W$.
Mich\`{e}le Vergne provided us with a counter-example.
She also informed us about the references \cite{Hsi} and
\cite{Duf-Kum-Ver}.
\end{note}

\begin{note} We prefer to use $u^a$ instead of $\Theta^a$ 
when the Lie group is a torus, 
and $u$ in the case of a circle.
\end{note}

\subsection{Equivariant Euler class} We will also need the notion of 
equivariant  Euler class \cite{Ati-Bot}.
Let $F$ be a connected closed oriented manifold, 
and $\pi: E \rightarrow F$ be a smooth complex 
vector bundle over $F$. 
Assume that there is an 
$S^1$-action on $E$ by bundle homomorphisms, 
which covers an $S^1$-action on $F$. 
Then one  can define \cite{Ati-Bot}
the equivariant Euler class $\epsilon(E) \in H^*_{S^1}(F)$,
which satisfies
$$\epsilon (E_1 \oplus E_2) = \epsilon(E_1) \epsilon(E_2)$$
for two $S^1$ bundles $E_1$ and $E_2$ over $F$. 
We will be concerned with the case when the action 
of $S^1$ on $F$ is trivial. In this case, 
$E$ has a decomposition as $S^1$ bundles 
$$E = L_1 \oplus L_2 \oplus \cdots L_r,$$
where each $L_j$ is a line bundle such that 
the action of $\exp (2\pi \sqrt{-1}t)$ on $L_j$ is 
multiplication by $\exp ( 2\pi \sqrt{-1}m_j t)$, 
for some weight $m_j \in {\mathbb Z}$. 
By formula $(8.8)$ in Atiyah-Bott \cite{Ati-Bot}, 
$$\epsilon(L_j) = m_j u + c_1(L_j).$$
Hence we have 
$$\epsilon(E) = \prod_{j=1}^r (m_j u  + c_1(L_j)).$$
Here $u$ is dual to an element $\xi$ in the Lie
algebra of $S^1$,
such that if $S^1$ is given an invariant metric 
in which $|\xi| =1$, then $\vol (S^1) = 1$.

\subsection{Borel model} Historically, equivariant cohomology was defined 
by Borel model (cf. Atiyah-Bott \cite{Ati-Bot}). 
Let $\pi: EG \rightarrow BG$ be a universal principal 
$G$-bundle (it  is unique up to homotopy).
The Borel constuction of a $G$-manifold $X$ is
$$X_G := (X \times EG)/G,$$
where $G$ acts on $X \times EG$ diagonally.
It can be shown that $H^*(X_G, {\Bbb R}) \cong H^*_G(X)$ 
(see e.g. Lawson \cite{Law}).
In particular, when $X$ is a point,
$X_G = BG$,
one has
$$H^*(BG, {\Bbb R}) \cong H^*_G (pt) \cong S(\fg^*)^G.$$
The homomorphism
$H^*(BG, {\Bbb R}) \rightarrow H^*(X_G, {\Bbb R})$ induced from
the map of $G$-spaces $X \rightarrow pt$ can be identified with 
the homomorphism 
$$S(\fg^*)^G \rightarrow H^*_G(X).$$

Since their relationships  with differential forms,
the Cartan model and Weil model became popular 
after the works of Berline-Vergne \cite{Ber-Ver1, Ber-Ver2},
Atiyah-Bott \cite{Ati-Bot}, Mathai-Quillen \cite{Mat-Qui}, etc.

\section{Degenerate Chern-Weil Theory}

 Our construction in this section is motivated by 
the equivariant Chern-Weil theory
for equivariant principal bundles 
due to Berline-Vergne \cite{Ber-Ver1}.
In the special case of $G = S^1$, the construction is used in 
Cao-Zhou \cite{Cao-Zho} to prove wall crossing formulas 
in Seiberg-Witten theory. 
(In fact, the original construction in \cite{Cao-Zho} was different and
more complicated,
we were led to the present version by consideration 
of generalization to the nonabelian case.)

Let $W$ be a compact $G$-manifold, 
possibly with boundary $\partial W$, 
such that the $G$-action on an open subset of $W$, 
which contains $\partial W$, is (locally) free. 
We will call $W$ a degenerate principal $G$-bundle. 
 Denote by $W^s$ the set of points in $W$
whose stabilizers have dimension $> 0$ and set $W^0 = W - W^s$.  
Let $f: W \rightarrow [0, 1]$ be a $G$-invariant smooth 
function on $W$
which  vanishes on a tubular neighborhood of $W^s$, 
and is identically $1$   outside a larger tubular neighborhood. 
Let $\omega$ be a connection of the principal bundle 
$W^0 \rightarrow W^0/G$.
We call 
$\omega_f = f \cdot \omega   \in \fg \otimes \Omega^1(W)$
a degenerate connection, and 
$$\Omega_f^G = d \omega_f + \frac{1}{2} [ \omega_f, \omega_f] - 
(-1 + f) \xi_a \Theta^a \in \fg \otimes \Omega^2_G(W)$$
the degenerate equivariant curvature of $\omega_f$.

\begin{lemma} We have $D \Omega_f^G = [\Omega_f^G, \omega_f]$.
\end{lemma}

\begin{proof} It suffices to prove it on $W^0$, 
on which we have  $\iota_a \omega = \xi_a$, $\iota_a df = 0$. 
Furthermore,
\begin{eqnarray*}
& & \iota_a d \omega = \iota_a (d \omega + 
	\frac{1}{2} [\omega, \omega] ) 
	- \frac{1}{2} \iota_a  [\omega, \omega] 
	=  - [\xi_a, \omega], \\
& & [[\omega_f,  \omega_f], \omega_f] = 
	f^3 [[\omega, \omega], \omega] = 0
\end{eqnarray*}
So on $W^0$, we have
\begin{eqnarray*}
D \Omega^G_f & = & 
	d ( d \omega_f + \frac{1}{2} [ \omega_f, \omega_f] - 
	(-1 + f) \xi_a \Theta^a) \\
& &     - \Theta^b \iota_b (d \omega_f 
	+ \frac{1}{2} [ \omega_f, \omega_f] - 
	(-1 + f) \xi_a \Theta^a) \\
& = & [d \omega_f, \omega_f] - df \xi_a \Theta^a 
 	- \Theta^a \iota_a (df \wedge \omega + f d \omega  
   	+ \frac{f^2}{2} [\omega, \omega]) \\
&  = & [d \omega_f + \frac{1}{2} [\omega_f, \omega_f], \omega_f] 
   - df \xi_a \Theta^a + \Theta^a df \xi_a 
   + (f - f^2)  \Theta^a [\xi_a, \omega] \\
&  = & [\Omega_f, \omega_f] 
   + [(1 -f) \xi_a \Theta^a, f \omega] = [\Omega_f^G, \omega_f].
\end{eqnarray*}

\end{proof}

In the remaining part of this section, we shall adopt the following

\begin{conventions}
If $i_1, \cdots, i_q$ are indices, then $(i_1 \cdots i_q)$ 
means symmetrizing on these indices, 
and $[i_1 \cdots i_q]$ means antisymmetrize on these indices. 
Furthermore, notation like $i_1 \cdots  |b| \cdots i_q$ means 
$b$ does not participate in the (anti-)symmetrization.
\end{conventions}

\begin{lemma} Let $F \in S^q(\fg^*)^G$,  then 
$$q F([\Omega_f^G, \omega_f], \Omega_f^G, \cdots, \Omega_f^G) 
= 0.$$
\end{lemma}

\begin{proof} As in the 
ordinary case, this is equivalent to the invariance of $F$.  
Let $F = a_{i_1  \cdots i_q} \Theta^{i_1} \cdots 
\Theta^{i_q}$, 
where $a_{i_1  \cdots i_q} = a_{(i_1  \cdots i_q)}$. 
Since $F$ is $G$-invariant, we have
\begin{eqnarray*}
0 & = & L_b F = q a_{i_1  \cdots i_q} (L_b \Theta^{i_1}) 
	\Theta^{i_2} \cdots  \Theta^{i_q} \\
& = & - q f^{i_1}_{bc} a_{i_1  \cdots i_q} \Theta^c 
	\Theta^{i_2} \cdots  \Theta^{i_q} \\
& = & - q f^{i_1}_{b(c} a_{|i_1|i_2  \cdots i_q)} \Theta^c 
	\Theta^{i_2} \cdots  \Theta^{i_q},	.
\end{eqnarray*}
where the last term is obtained after symmetrizing the indices 
$c, i_2,\cdots, i_q$. Therefore, $q f^{i_1}_{b(c} a_{|i_1|i_2  \cdots i_q)} = 
0$.  
Hence we have
\begin{eqnarray*}
& & q F([\Omega_f^G, \omega_f], \Omega^G_f, \cdots 
	\Omega_f^G) \\
& = & - \omega_f^b q f^{i_1}_{bc} a_{i_1  \cdots i_q} 
	(\Omega_f^G)^c \wedge (\Omega_f^G)^{i_2} 
	\wedge \cdots  \wedge (\Omega_f^G)^{i_q} \\
& = & - \omega_f^b q f^{i_1}_{b(c} a_{|i_1|i_2  \cdots i_q)} 
	(\Omega_f^G)^c \wedge (\Omega_f^G)^{i_2} 
	\wedge \cdots  \wedge (\Omega_f^G)^{i_q}
	= 0.
\end{eqnarray*}

\end{proof}

\begin{lemma} Let $F \in S^q(\fg^*)^G$, then 
$F(\Omega_f^G, \cdots, \Omega_f^G) \in (S(\fg^*) \otimes \Omega(W))^G$.
Furthermore,
$$D F(\Omega_f^G, \cdots, \Omega_f^G) = 0.$$
\end{lemma}

\begin{proof}
The first statement is obvious. For the second,
\begin{eqnarray*}
& & D F(\Omega_f^G, \cdots, \Omega_f^G) 
    = q F (D\Omega_f^G, \Omega_f^G, \cdots, \Omega_f^G) \\
& = & q F([\Omega_f^G, \omega_f], \Omega_f^G, 
	\cdots, \Omega_f^G) = 0.
\end{eqnarray*}

\end{proof}

As a corollary to Lemma 2.3, we have

\begin{theorem}  Let $W$ be a degenerate principal $G$-bundle. Given
a degenerate connection $\omega_f$ with equivariant degenerate 
curvature $\Omega_f^G$, there is a homomorphism (called degenerate 
Chern-Weil homomorphism)
$$cw_f: S(\fg^*)^G \rightarrow H^*_G(W),$$
which is induced from the homomorphism 
$$CW_f: S(\fg^*) \rightarrow S(\fg^*) \otimes \Omega(W)$$
given by $F \in S(\fg^*) \mapsto F(\Omega_f^G)$. 
\end{theorem}

Similar to Lemma 2.2,  one can prove the following

\begin{lemma}
Let $F \in S^q(\fg^*)^G$,  then 
$$q (q-1) F(\alpha, [\Omega_f^G, \omega_f], \Omega_f^G, \cdots, \Omega_f^G) = 
q F([\alpha, \omega_f], \Omega_f^G, \cdots, \Omega_f^G)$$
for $\alpha \in \Omega^1(W) \otimes \fg$.
\end{lemma}

\begin{theorem} The degenerate Chern-Weil homomorphism in Theorem 2.1
does not depend on the choice of the connection $\omega$ on $W^0$ or
the cut-off function $f$.
\end{theorem}

\begin{proof}
Let $\omega^0_f$ and $\omega^1_f$ be two connections on $W$, with degenerate
equivariant curvatures 
$(\Omega_f^G)^0$ and  $(\Omega_f^G)^1$ respectively,
  and consider
$$\widetilde{\omega}_f = (1 - t) \omega_f^0 + t \omega_f^1.$$
Then $\widetilde{\omega}_f$ is a degenerate connection on $W \times I$, $I = [0, 1]$. 
Denote by $\widetilde{\Omega}_f^G$
the degenerate equivariant 
curvature of $\widetilde{\omega}_f$, $\pi: W \times I \rightarrow W$ the 
projection, and let
$$\int_{\pi}: W(\fg) \otimes \Omega(W \times I) \rightarrow 
W(\fg) \otimes \Omega(W)$$ 
be defined by
$$\int_{\pi} \alpha(t) +  dt \wedge \beta(t)= \int_0^1 \beta(t) dt,$$
where  $\alpha(t)$ and $\beta(t)$ are families of 
equivariant differential forms on
$W$ depending smoothly on $t$. For any $F \in S^q(\fg^*)^G$, define the
degenerate transgression operator 
$$T_{(\omega_f^0, \omega_f^1)} F = \int_{\pi} F(\widetilde{\Omega}_f^G).$$
Then one can check that
$$D T_{(\omega_f^0, \omega_f^1)} F = F((\Omega_f^G)^1) - F((\Omega_f^G)^0).
$$
Indeed, if we let $\delta = \omega^1_f - \omega^0_f$ and
$(\Omega_f^G)^t$ be the degenerate equivariant curvature of
$\omega_f^t = \omega^0_f + t \delta$, then
$$(\Omega_f^{\mbG})^t = (\Omega^{\mbG}_f)^0 + td\delta + 
t[\omega^0_f, \delta] + \frac{t^2}{2}[\delta, \delta] + dt \wedge \delta.$$
So we have
\begin{eqnarray*} 
& \frac{d}{dt} (\Omega_f^G)^t = 
	d\delta + [\omega_f^t, \delta] = 
	D \delta + [\omega_f^t, \delta], \\
& T_{(\omega_f^0, \omega_f^1)} F 
= \int_0^1 q F(\delta, (\Omega_f^G)^t \cdots, (\Omega_f^G)^t) dt.
\end{eqnarray*}
Hence,
\begin{align*}
  &   D T_{(\omega_f^0, \omega_f^1)} F  && \\
 = &  \int_0^1 (q F(D \delta, (\Omega_f^G)^t,
	\cdots, (\Omega_f^G)^t)  dt &&\\
   & +	q(q-1) F(\delta, D(\Omega_f^G)^t, (\Omega_f^G)^t,
	\cdots, (\Omega_f^G)^t) ) dt&&\\
 = & \int_0^1 (q F(D \delta, (\Omega_f^G)^t,
	\cdots, (\Omega_f^G)^t) &&\\
  &  + q(q-1) F(\delta, [(\Omega_f^G)^t, \omega_f^t],
 	(\Omega_f^G)^t, \cdots, (\Omega_f^G)^t) ) dt 
 	&& \text{(by Lemma 2.1)} \\
 = & \int_0^1 q F(D \delta + [\omega_f^t, \delta], (\Omega_f^G)^t,
	\cdots, (\Omega_f^G)^t) dt  && \text{(by Lemma 2.4)} \\
 = & \int_0^1 q p(\frac{d}{dt}(\Omega_f^G)^t, (\Omega_f^G)^t,
	\cdots, (\Omega_f^G)^t) dt &&\\
 = & \int_0^1 \frac{d}{dt} F((\Omega_f^G)^t, (\Omega_f^G)^t,
	\cdots, (\Omega_f^G)^t) dt &&\\
 = &  F((\Omega_f^G)^1, \cdots, (\Omega_f^G)^1) - 
	F((\Omega_f^G)^0, \cdots, (\Omega_f^G)^0).&& 
\end{align*}
Therefore, the degenerate Chern-Weil homomorphism is independent of the 
choice of $\omega_f$. Similarly, if $f^0$ and $f^1$
are two cut-off functions used to carry out the construction, then
on $W \times I$, setting $\bar{f} = (1 - t)f^0 + t f^1$,
 $\overline{\omega}_{\bar{f}} = {\bar{f}} \omega$, and
$$\overline{\Omega}_{\bar{f}}^G = d \overline{\omega}_{\bar{f}}+
\frac{1}{2}[\overline{\omega}_{\bar{f}}, \overline{\omega}_{\bar{f}}] 
- (-1 + \bar{f}) \xi_a\Theta^a,$$ 
we can define a similar transgression opertaor:
$$T_{(f^0, f^1)} F
	= \int_{\pi} F(\overline{\Omega}_{\bar{f}}^G).$$
Then the same proof as above shows that
$$D T_{(f^0, f^1)} F = F(\Omega_{f^1}^G) - F(\Omega_{f^0}^G).$$
Hence the degenerate Chern-Weil homomorphism is also 
independent of the choice of $f$.
\end{proof}

 Now consider the homomorphism
$CW_f^w: W(\fg) \otimes \Omega(W) \rightarrow 
S(\fg^*) \otimes \Omega(W)$
defined by extending the Chern-Weil construction
$\theta^a \mapsto  \omega_f^a$, 
$ \Theta^a  \mapsto (\Omega^G_f)^a$ as a $\Omega(W)$-module map.
 Define 
$r_f^G: S(\fg^*) \otimes \Omega(W) \rightarrow 
S(\fg^*) \otimes \Omega(W)$
by
$$r_f^G(\alpha) = CW_f^w(\Psi(\alpha)).$$
Let $U$ be any open set on which $f = 1$, then on $U$, 
we have $\omega_f = \omega$, $\Omega_f^G = \Omega$.
Therefore $r_f^G = r^G$ on $U$. 
So $r_f^G$ is a generalization of $r^G$. 
It is easy to see that $r_f^G$ maps 
$(S(\fg^*) \otimes \Omega(W))^G$ to itself.  

\begin{theorem} The homomorphism 
$r_f^G: (S(\fg^*) \otimes \Omega(W))^G \rightarrow
 (S(\fg^*) \otimes \Omega(W))^G$
satisfies $D r_f^G = r_f^G D$. 
Hence it induces a homomorphism of cohomologies: 
$$(r_f^G)_*: H^*_G(W) \rightarrow H^*_G(W).$$
\end{theorem}

\begin{theorem}
The homomorphism $(r_f^G)_*$ in Theorem 2.3 
does not depend on the choice of
the connection $\omega$ on $W^0$ 
or the choice of the cut-off function $f$.
\end{theorem}

An important observation, pointed out to us by
Mich\`{e}le Vergne, is that, 
if one takes $f \equiv 0$, then 
$\Omega_f^G = \xi_a\Theta^a$. 
Therefore, we have

\begin{theorem}
The homomorphism $(r^G_f)_*$ on cohomology is the identity map.
\end{theorem}

Let $i^c$ denote the inclusion $S(\fg^*) \rightarrow 
S(\fg^*) \otimes \Omega(W)$. 
It is obvious that $CW_f = r_f^G \circ i^c$. 
This equality reveals that Theorem 2.1 and Theorem 2.2 
are special cases of Theorem 2.3 and Theorem 2.4 respectively. 
Since $cw_f = (r_f^G)_* \circ i^c_*$, we have the following

\begin{theorem} We have $cw_f = i^c_*$.
\end{theorem}

We now present the proof of Theorems 2.3 and 2.4.
Theorem 2.4 is an easy consequence of Theorem 2.3
by a construction of transgression homomorphism.
Our proof of Theorem 2.3 relies on calculations by brute force.
It will be nice to  find a more conceptual proof.
To begin with, 
we have the following lemma which plays a similar role in the proof 
of Theorem 2.3 as  Lemma 2.4 in the proof of Theorem 2.1.

\begin{lemma} Let 
$\alpha = \Theta^{i_1} \cdots \Theta^{i_q} \alpha_{i_1\cdots i_q}
\in (S(\fg^*) \otimes \Omega(X))^G$, 
$\alpha_{i_1\cdots i_q} = \alpha_{(i_1\cdots i_q)}$, then we have
\begin{equation}
{\cal L}_b \alpha_{i_1 \cdots i_q} =  
      d \iota_b \alpha_{i_1 \cdots i_q} + \iota_b d \alpha_{i_1 \cdots i_q} 
 =  q f^p_{b(i_1} \alpha_{|p|i_2 \cdots  i_q)}.
\end{equation}
Furthermore, if  $\alpha$ is $D$-closed then we have
\begin{equation}
d \alpha_{i_1 \cdots i_q} = \iota_{(i_1} \alpha_{i_2 \cdots i_q)}.
\end{equation}
\end{lemma}

\begin{proof}
By $G$-invariance of $\alpha$,
\begin{eqnarray*}
0 & = & (L_b \otimes 1 + 1 \otimes {\cal L}_b) 
	\sum \Theta^{i_1} \cdots \Theta^{i_q} 
	\alpha_{i_1 \cdots i_q} \\
& = & \sum \Theta^{i_1} \cdots \Theta^{i_q} {\cal L}_b 
 	\alpha_{i_1 \cdots i_q}
	+ \sum L_b(\Theta^{i_1} \cdots \Theta^{i_q}) 
	\alpha_{i_1 \cdots i_q} \\
& = & \sum \Theta^{i_1} \cdots \Theta^{i_q} {\cal L}_b 
	\alpha_{i_1 \cdots i_q}
  + \sum q L_b(\Theta^{i_1}) \Theta^{i_2}  \cdots \Theta^{i_q} 
  	\alpha_{i_1 \cdots i_q} \\
& = & \sum \Theta^{i_1} \cdots \Theta^{i_q} {\cal L}_b 
	\alpha_{i_1 \cdots i_q}
  - \sum q f_{bp}^{i_1}\Theta^p \Theta^{i_2}  \cdots \Theta^{i_q} 
  	\alpha_{i_1 \cdots i_q} \\
& = & \sum \Theta^{i_1} \cdots \Theta^{i_q} 
	({\cal L}_b \alpha_{i_1 \cdots i_q} - q f^p_{bi_1} 
	\alpha_{pi_2 \cdots  i_q}).
\end{eqnarray*}
The last equality is obtained by interchanging $p$ with $i_1$. 
This proves $(4)$. Similarly, from $D \alpha = 0$, we get
\begin{eqnarray*}
0 & = & \sum \Theta^{i_1} \cdots \Theta^{i_q} 
	d \alpha_{i_1 \cdots i_q}
        - \sum \Theta^{i_1} \cdots \Theta^{i_q} \Theta^b 
        i_b \alpha_{i_1 \cdots i_q} \\
  & = & \sum \Theta^{i_1} \cdots \Theta^{i_q}
  	(d \alpha_{i_1 \cdots i_q} - 
  	\iota_{(i_1} \alpha_{i_2 \cdots i_q)}).
\end{eqnarray*} 
This proves $(5)$.
\end{proof}

\begin{proof}[Proof of Theorem 2.3] 
Let $\alpha =  \Theta^{I} \alpha^q_I
	= \sum_q \Theta^{i_1} \cdots \Theta^{i_q} 
	\alpha_{i_1 \cdots i_q}$,
then using the summation convention, we have
\begin{eqnarray*}
& & \Psi (\alpha) = \prod_{a = 1}^k 
	(1 - \theta^a \otimes \iota_a) \alpha \\
& = & \sum_{j = 0}^k (-1)^{j(j+1)/2} 
	\sum_{a_1 < \cdots < a_j} 
	\theta^{a_1} \cdots \theta^{a_j} 
 	\Theta^{I} \iota_{a_1} \cdots \iota_{a_j} \alpha_I  \\
&  = & \sum_{j = 0}^k \frac{(-1)^{j(j+1)/2} }{j!} 
	\theta^{a_1} \cdots \theta^{a_j} 
 	\Theta^{I} \iota_{a_1} \cdots \iota_{a_j} \alpha_I.
\end{eqnarray*}
Applying the degenerate Chern-Weil construction $CW_f$, we get
\begin{eqnarray*}
r_f^G (\alpha) =
\sum_{j = 0}^k  \frac{(-1)^{j(j+1)/2}}{j!} 
	\omega_f^{a_1} \wedge  \cdots 
	\wedge \omega_f^{a_j} \wedge
 	(\Omega_f^G)^{I} \wedge \iota_{a_1} \cdots \iota_{a_j} 
 	\alpha_{\mbox{\tiny $I$}}.
\end{eqnarray*}
Taking $D$ on both sides, we see that $D r_f^G (\alpha)$ is equal to:
\begin{align} 
  & \sum_{j = 0}^k  \frac{(-1)^{j(j+1)/2}}{j!} 
 	D (\omega_f^{a_1} \wedge  \cdots 
	\wedge \omega_f^{a_j}) \wedge
 	(\Omega_f^G)^{I} \wedge \iota_{a_1} 
 	\cdots \iota_{a_j} \alpha_I \\
 + & \sum_{j = 0}^k  \frac{(-1)^{j(j+3)/2}}{j!} 
 	\omega_f^{a_1} \wedge  \cdots 
	\wedge \omega_f^{a_j} \wedge
 	D (\Omega_f^G)^{I} \wedge \iota_{a_1} 
 	\cdots \iota_{a_j} \alpha_I \\
 + & \sum_{j = 0}^k  \frac{(-1)^{j(j+3)/2}}{j!} 
 	\omega_f^{a_1} \wedge  \cdots 
	\wedge \omega_f^{a_j} \wedge
 	(\Omega_f^G)^{I} \wedge D (\iota_{a_1} 
 	\cdots \iota_{a_j} \alpha_I). 
\end{align}
We will examine each of the above terms separately. 
To start with, recall that
\begin{eqnarray*}
D \omega_f^a & = & d \omega_f^a - f \Theta^a \\
& = & (\Omega^G_f)^a 
	- \frac{1}{2} f^{a}_{bc} 
	\omega_f^b \wedge  \omega_f^c - \Theta^a, \\
D (\Omega^G_f)^a & = & 
	f^a_{bc} (\Omega_f^G)^b \wedge \omega_f^c
	=- f^a_{bc} \omega_f^b (\Omega_f^G)^c.
\end{eqnarray*}
Then $(6)$ can be written  as

\begin{align}
   & \sum_{j = 0}^k \frac{(-1)^{j(j+1)/2}}{j!} 
   	D (\omega_f^{a_1} \wedge  \cdots \wedge \omega_f^{a_j}) 
   	\wedge (\Omega_f^G)^{I} \wedge 
   	\iota_{a_1} \cdots \iota_{a_j} \alpha_I \nonumber\\	
 = & \sum_{j = 0}^k  \frac{(-1)^{j(j+1)/2}}{j!} 
 	j D \omega_f^{a_1} \wedge \omega_f^{a_2} \wedge  \cdots 
 	\wedge \omega_f^{a_j} \wedge
 	(\Omega_f^G)^{I} \wedge 
 	\iota_{a_1} \cdots \iota_{a_j} \alpha_I \nonumber \\
 = & \sum_{j = 1}^k  \frac{(-1)^{j(j+1)/2}}{(j - 1)!} 
	((\Omega^G_f)^{a_1}
	- \frac{1}{2} f^{a_1}_{bc} \omega_f^b \wedge  \omega_f^c 
	- \Theta^{a_1})
	\wedge \omega_f^{a_2} \wedge  \cdots \wedge \omega_f^{a_j} 
	\nonumber \\
 & 	\wedge
 	(\Omega_f^G)^{I} \wedge \iota_{a_1} \cdots \iota_{a_j} 
 	\alpha_I \nonumber \\
 = &  \sum_{j = 1}^k \frac{(-1)^{j(j+1)/2}}{(j-1)!} 
 	\omega_f^{a_2} \wedge  \cdots 
	\wedge \omega_f^{a_j} \wedge
 	(\Omega_f^G)^{I} \wedge (\Omega_f^G)^{a_1} 
 	\wedge \iota_{a_1} \cdots \iota_{a_j} \alpha_I  
 	\tag{6a} \\
 - & \sum_{j = 1}^k \frac{(-1)^{j(j+1)/2}}{2(j-1)!} 
 	f^{a_1}_{bc} \omega_f^b \wedge 
	\omega_f^c \wedge \omega_f^{a_2} \wedge  \cdots 
	\wedge \omega_f^{a_j} \wedge
 	(\Omega_f^G)^{I} \wedge  \iota_{a_1}  \cdots \iota_{a_j} 
 	\alpha_I  \tag{6b} \\ 	
 - & \sum_{j = 1}^k \frac{(-1)^{j(j+1)/2}}{(j-1)!} 
 	\omega_f^{a_2} \wedge  \cdots 
	\wedge \omega_f^{a_j} \wedge
 	(\Omega_f^G)^{I} \wedge \Theta^{a_1} 
 	\iota_{a_1} (\iota_{a_2} \cdots \iota_{a_j} \alpha_I) 	\tag{6c}
\end{align}	 	
Similarly, we rewrite $(8)$ as 
\begin{align}
& \sum_{j = 0}^k  \frac{(-1)^{j(j+3)/2}}{j!} 
   	\omega_f^{a_1} \wedge  \cdots 
	\wedge \omega_f^{a_j} \wedge (\Omega_f^G)^{I} \wedge 
	D (\iota_{a_1} \cdots \iota_{a_j} \alpha_I) \nonumber \\
= & \sum_{j = 0}^k  \frac{(-1)^{j(j+3)/2}}{j!} 
 	\omega_f^{a_1} \wedge  \cdots 
	\wedge \omega_f^{a_j} \wedge
 	(\Omega_f^G)^{I} \wedge 
 	d (\iota_{a_1} \cdots \iota_{a_j} \alpha_I) \tag{8a} \\
- & \sum_{j = 0}^k  \frac{(-1)^{j(j+3)/2}}{j!} 
	\omega_f^{a_1} \wedge  \cdots 
	\wedge \omega_f^{a_j} \wedge
 	(\Omega_f^G)^{I} \wedge \Theta^b 
 	\iota_b (\iota_{a_1} \cdots \iota_{a_j} \alpha_I) 
 	\tag{8b}
\end{align}
By a renaming of the indices, it is
easy to see that $(6c)$ and $(8b)$ together yield
$$\sum_{b = 1}^k -(-1)^{k(k+3)/2} \omega_f^1 \wedge  \cdots 
	\wedge \omega_f^k \wedge
 	(\Omega_f^G)^{I} \wedge \Theta^b 
 	\iota_b (\iota_1 \cdots \iota_k \alpha_I) = 0.$$
Now $(7)$ can be written as summation for $j = 0$ to $k$ of 
$\frac{(-1)^{j(j+3)/2}}{j!} \omega_f^{a_1} \wedge  \cdots 
	\wedge \omega_f^{a_j}$ wedge the following terms
	
\begin{eqnarray*}
&   & D (\Omega_f^G)^{I} \wedge \iota_{a_1} \cdots \iota_{a_j} \alpha_I 
 	\nonumber \\
& = &  D ((\Omega_f^G)^{i_1} \wedge \cdots \wedge (\Omega_f^G)^{i_q} ) \wedge 
 	\iota_{a_1} \cdots \iota_{a_j} \alpha_{i_1 \cdots i_q} \\
& = & q D (\Omega_f^G)^{i_1} \wedge (\Omega_f^G)^{i_2} \wedge \cdots
 	\wedge (\Omega_f^G)^{i_q} \wedge 
 	\iota_{a_1} \cdots \iota_{a_j} \alpha_{i_1 \cdots i_q} \\
& = & - q f^{i_1}_{bc} \omega_f^b \wedge (\Omega_f^G)^c \wedge 
 	(\Omega_f^G)^{i_2} \wedge \cdots \wedge (\Omega_f^G)^{i_q} \wedge
 	\iota_{a_1} \cdots \iota_{a_j} \alpha_{i_1 \cdots i_q} \\ 
& = & - \omega_f^b \wedge (\Omega_f^G)^c \wedge 
	(\Omega_f^G)^{i_2} \wedge \cdots \wedge 
 	(\Omega_f^G)^{i_q} \wedge
 	\iota_{a_1} \cdots \iota_{a_j} 
 	(q f^{i_1}_{bc} \alpha_{i_1 \cdots i_q}). \\		
& = & - \omega_f^b \wedge (\Omega_f^G)^c \wedge 
	(\Omega_f^G)^{i_2} \wedge \cdots \wedge 
 	(\Omega_f^G)^{i_q} \wedge
 	\iota_{a_1} \cdots \iota_{a_j} 
 	(q f^{i_1}_{b(c} \alpha_{|i_1|i_2 \cdots i_q)}) \\
& = & - \omega_f^b \wedge (\Omega_f^G)^c \wedge 
	(\Omega_f^G)^{i_2} \wedge \cdots \wedge 
 	(\Omega_f^G)^{i_q} \wedge
 	\iota_{a_1} \cdots \iota_{a_j} 
 	{\cal L}_b  \alpha_{c i_2 \cdots i_q}	\\
& = & - \omega_f^b \wedge (\Omega_f^G)^I \wedge 
	\iota_{a_1} \cdots \iota_{a_j} {\cal L}_b  \alpha_I.
\end{eqnarray*}
We have used $(4)$ in the second to last equality. 
Renaming $b$ by $a_1$, $a_l$ by $a_{l+1}$,
 one sees that $(7)$ is equal to 
$$ 
- \sum_{j = 0}^k \frac{(-1)^{j(j+1)/2}}{j!}
	\omega_f^{a_1} \wedge  \cdots \wedge \omega_f^{a_{j+1}} \wedge
	(\Omega_f^G)^I \wedge  \iota_{a_2} \cdots \iota_{a_{j+1}} 
 	{\cal L}_{a_1}  \alpha_I.
$$
The contribution from $j = k$ is clearly zero, so by changing  $j$ to
$j-1$,  $(7)$ is equal to 
\begin{equation}
- \sum_{j=1}^k \frac{(-1)^{j(j-1)/2}}{(j-1)!}
	\omega_f^{a_1} \wedge  \cdots \wedge \omega_f^{a_j} \wedge
	(\Omega_f^G)^I \wedge \iota_{a_2} \cdots \iota_{a_j} 
 	{\cal L}_{a_1} \alpha_I.   \tag{$7'$}
\end{equation}
Notice that in $(6a)$, it won't change the result if we
take the summation for $j = 1$ to $k+1$. Change the index $a_1$ to $b$, 
$a_l$ to $a_{l-1}$ for $l > 1$, one sees that $(6a)$ is equal to
\begin{equation}
- \sum_{j = 0}^{k} \frac{(-1)^{(j+1)j/2}}{j!} \omega_f^{a_1} \wedge  \cdots 
	\wedge \omega_f^{a_j} \wedge
 	(\Omega_f^G)^{I} \wedge  (\Omega_f^G)^b \wedge 
 	\iota_{a_1} \cdots \iota_{a_j} \iota_b \alpha_I.   
 	\tag{$6a'$}
\end{equation}	
Similarly,  $(6b)$ is equal to
\begin{equation}
-  \sum_{j = 2}^k \frac{(-1)^{j(j-1)/2}}{2(j-2)!}  
	\omega_f^{a_1} \wedge  \wedge  \cdots \wedge \omega_f^{a_j} \wedge
 	(\Omega_f^G)^{I} \wedge  
 	f^{c}_{a_1a_2}  \iota_{a_3}  \cdots \iota_{a_j} \iota_c \alpha_I.
\tag{$6b'$} \end{equation}
To summarize, we have $D r_f^G \alpha =  (6a') + (6b') + (7') +(8a)$.
On the other hand, 
\begin{align}
 & r_f^G D \alpha \nonumber \\
 = & \sum_{j = 0}^k  \frac{(-1)^{j(j+1)/2}}{j!} \omega_f^{a_1} \wedge  \cdots 
	\wedge \omega_f^{a_j} \wedge
 	(\Omega_f^G)^{I} \wedge \iota_{a_1} \cdots \iota_{a_j} 
 	d \alpha_{\mbox{\tiny $I$}} 
 	\tag{$9a$} \\
 - & \sum_{j = 0}^k  \frac{(-1)^{j(j+1)/2}}{j!} \omega_f^{a_1} \wedge  \cdots 
	\wedge \omega_f^{a_j} \wedge
 	(\Omega_f^G)^{I} \wedge (\Omega_f^G)^b \wedge
 	\iota_{a_1} \cdots \iota_{a_j} 
 	\iota_b \alpha_{\mbox{\tiny $I$}} \tag{$9b$}
\end{align} 
Since $(6a')$ cancels $(9b)$,  
$Dr_f^G \alpha - r_f^G D \alpha = (8a) - (9a) + (7') + (6b')$.
It is the summation for $j = 0$ to $k$ of
$\frac{(-1)^{j(j-1)/2}}{j!} \omega_f^{a_1} \wedge  \cdots 
	\wedge \omega_f^{a_j} \wedge
 	(\Omega_f^G)^{I} $ wedge 
$$(d \iota_{[a_1} \cdots \iota_{a_j]} 
	- (-1)^j \iota_{[a_1} \cdots \iota_{a_j]} d
 	- j \iota_{[a_2} \cdots \iota_{a_j} {\cal L}_{a_1]} - \frac{1}{2} j 
(j-1) 
 f^{c}_{[a_1a_2}  \iota_{a_3}  \cdots \iota_{a_j]} \iota_c) \alpha_I,$$
which is easily shown to be zero by
repeatedly using $d \iota_a = {\cal L}_a -  \iota_a d$
and ${\cal L}_a \iota_b = \iota_b {\cal L}_a - f^{c}_{ab}  \iota_c$. 
\end{proof}

\begin{lemma}
Let $D_W: \Omega_G(W) \rightarrow \Omega_G(W)$ and 
$D_{W \times I}: \Omega_G(W \times I) 
\rightarrow \Omega_G(W \times I)$
be the Cartan differentials on $W$ and $W \times I$ 
respectively. 
If $\widetilde{\alpha} \in \Omega_G(W \times I)$ can be 
written as
$$\widetilde{\alpha} = 
\alpha(t) +  dt \wedge \beta(t),$$
where  $\alpha(t)$ and $\beta(t)$ are families of 
equivariant differential forms on
$W$ which depend smoothly on $t$, then
$$\int_{\pi} D_{W \times I} (\widetilde{\alpha}) +
D_W \int_{\pi} \widetilde{\alpha} = \alpha(1) - \alpha(0).
$$
\end{lemma}

\begin{proof} Notice that 
$D_{W \times I} = D_W + dt \wedge \frac{\partial}{\partial t}$. Then we have
\begin{eqnarray*}
&   & \int_{\pi} D_{W \times I} (\widetilde{\alpha})
 	=\int_{\pi} D_{W \times I} (\alpha(t) +  dt \wedge \beta(t)) \\
& = & \int_{\pi} D_W \alpha(t) + 
	dt \wedge \frac{\partial}{\partial t} \alpha (t)  
	- dt \wedge D_W \beta (t) \\
& = & \int_0^1 [\frac{d}{dt} \alpha (t) - D_W \beta(t) ] dt
\end{eqnarray*}
On the other hand,
$$D_W \int_{\pi} \widetilde{\alpha} 
	= D_W \int_0^1 \beta(t) dt \\
	= \int_0^1 [D_W \beta (t)] dt.$$
So we have
$$ \int_{\pi} D_{W \times I} (\widetilde{\alpha}) +
	D_W \int_{\pi} \widetilde{\alpha} 
	 =  \int_0^1 \frac{d}{dt} \alpha (t) dt
	 = \alpha(1) - \alpha(0).
$$
\end{proof}

\begin{proof}[Proof of Theorem 2.4] 
For any $\alpha \in \Omega_G(W)$, we have 
$$D_{W \times I} (\pi^*\alpha) = \pi^* (D_W \alpha).$$
For two connections 
$\omega^0$ and $\omega^1$ on $W^0$,  
use the notations 
in the proofs of Theorem 2.2 and Theorem 2.3, 
we define a degenerate transgression operator
\begin{eqnarray*}
&   & T_{(\omega^0, \omega^1)}\alpha \\
& = & \int_{\pi} \sum_{j = 0}^k  \frac{(-1)^{j(j+1)/2}}{j!} 
	\widetilde{\omega}_f^{a_1} \wedge  \cdots 
	\wedge \widetilde{\omega}_f^{a_j} \wedge
 	(\widetilde{\Omega}_f^G)^{I} \wedge 
 	\pi^* (\iota_{a_1} \cdots \iota_{a_j} \alpha_I) \\
& = & \int_{\pi} \widetilde{r}_f^G (\pi^* \alpha)
\end{eqnarray*}
where $\widetilde{r}_f^G: \Omega_G(W \times I) \rightarrow \Omega_G(W \times I)$
is obtained by the degenerate Chern-Weil construction for 
$\widetilde{\omega}_f$ and $\widetilde{\Omega}_f^G$. Then by Lemma 2.6 and
Theorem 2.3 for $W \times I$,
we have
\begin{align*}
   & D_W T_{(\omega^0, \omega^1)} \alpha 
   	+ T_{(\omega^0, \omega^1)} D_W \alpha  && \\
= & D_W \int_{\pi} \widetilde{r}_f^G (\pi^* \alpha)
 	+ \int_{\pi} \widetilde{r}_f^G (\pi^* D_W \alpha)  && \\
= & D_W \int_{\pi} \widetilde{r}_f^G (\pi^* \alpha)
 	+ \int_{\pi} \widetilde{r}_f^G D_{W \times I} 
 	(\pi^*  \alpha)  && \\	
= & D_W \int_{\pi} \widetilde{r}_f^G (\pi^* \alpha)
 	+ \int_{\pi} D_{W \times I} (\widetilde{r}_f^G 
 	(\pi^*  \alpha))  && \\ 
= & (r_f^G)^1 (\alpha) - (r_f^G)^0 (\alpha) 
	&& \text{(by Lemma 2.6)}
\end{align*}
So $(r_f^G)_*$ is independent of the choice of $\omega$. 
A similar construction shows that it is also independent of the choice of $f$.
\end{proof}

\section{Two nonabelian localization formulas}

In this section, we prove two nonabelian localization formulas of Kalkman
type for $G = SU(2)$ and $SO(3)$. 
First, let us recall Kalkman's localization formula for circle action 
\cite{Kal2} stated as follows:

\begin{proposition}
Let $W$ be an $S^1$-manifold with an invariant boundary 
$\partial W$, 
such that the $S^1$-action on $\partial W$ is locally free 
and effective.
 Suppose that
$F = \{ P_k \}$ is a decomposition of the fixed point set 
 into connected components.
Denote by $\nu_k$ the normal bundle of $P_k$
in $W$, and $\epsilon(\nu_k)$ the equivariant Euler class
of $\nu_k$. Then for any homogeneous $D_{S^1}$-closed form
 $\alpha$  on $W$ of total degree $dim (W) - 2$, we have
$$\int_{\partial W/ S^1} r^{S^1}(\alpha) = 
\sum_k \int_{P_k} \frac{\alpha u}{\epsilon (\nu_k)}.$$
\end{proposition}

In the above formula, we have used the normalization such that
$S^1$ has volume $1$.
In our earlier work \cite{Cao-Zho}, we have applied this formula
to obtain wall crossing formulas in Seiberg-Witten theory
due to Li-Liu \cite{Li-Liu} and Okonek-Teleman \cite{Oko-Tel}.
We now state our first nonabelian generalization:

\begin{theorem} \label{formula1}
Assume that $G = SU(2)$ or $SO(3)$ 
acts on a compact manifold $W$ with boundary $\partial W$, 
such that the $G$-action on $\partial W$ is locally free and effective. 
Let $T \subset G$ be a circle subgroup,
with fixed point set $F = \{P_k\}$. 
Suppose that $G$ is given a bi-invariant metric such 
that $\vol (T) = 1$. Then for any homogeneous $D_{\mbox{\tiny $G$}}$-closed 
form $\alpha$ of total degree $\dim (W) - 4$, we have
\begin{equation}
\int_{\partial W/G} r^{\mbox{\tiny $G$}}(\alpha) = - \frac{1}{c(G)} 
\sum_k \int_{P_k} \frac{p(\alpha) u^2}{\epsilon(\nu_k)},
\end{equation}
where $c(SU(2)) = 1$, $c(SO(3)) = 2$ and $p: \Omega_G(W) 
\rightarrow \Omega_T(W)$ is defined in $(3)$.
\end{theorem}

We will  apply Theorem \ref{formula1} to wall crossing in symplectic
reduction in a forthcoming paper \cite{Cao-Zho2}.
Before we embarking on a proof of Theorem \ref{formula1},
let us explain why one can expect localization formula
on manifolds.
The general situation is as follows.
Let $W$ be a compact $n$-dimensional manifold with boundary,
$G$ a $m$-dimsnional compact connected Lie group acting on $W$,
such that he acton on $\partial W$ is free,
$\alpha$ an equivariant closed form of degree $n-m-1$ on $W$.  
The problem is to compute 
$\int_{\partial W/G} r^G(\alpha|_{\partial W})$.
Suppose that $\fg$ is given a $G$-invariant Euclidean metric,
$\{\xi_1, \cdots, \xi_m\}$ an orthonrmal basis
with structure constants $c^j_{kl}$'s,
$\omega_f = \sum_a \omega_f^a \xi_a$ a degenerate connection on $W$.
Then by Stokes theorem and $Dr^G_f(\alpha) = 0$,
we  have
\begin{eqnarray*}
&   & \int_{\partial W/G} r^G(\alpha|_{\partial W})
= \frac{1}{\vol(G)} \int_{\partial W} 
\omega^1_f \wedge \cdots \wedge \omega^m_f \wedge r^G_f(\alpha) \\
& = & \frac{1}{\vol(G)} \int_W 
d(\omega^1_f \wedge \cdots \wedge \omega^m_f \wedge r^G_f(\alpha))
=\frac{1}{\vol(G)} \int_W 
D(\omega^1_f \wedge \cdots \wedge \omega^m_f \wedge r^G_f(\alpha)) \\
& = & \frac{1}{\vol(G)} \int_W \sum_{j=1}^m
(-1)^{j-1} (D\omega_f^j) \wedge 
\omega^1_f \wedge \cdots \wedge \widehat{\omega^j_f} \wedge \cdots
\wedge \omega_f^m \wedge r^G_f(\alpha) \\
& = & \frac{1}{\vol(G)} \int_W \sum_{j=1}^m
(-1)^{j-1} 
(\Omega_f^j - \frac{1}{2}c^j_{kl} \omega_f^k\wedge \omega_f^l) 
\wedge 
\omega^1_f \wedge \cdots \wedge \widehat{\omega^j_f} \wedge \cdots
\wedge \omega_f^m \wedge r^G_f(\alpha) \\
& = & \frac{1}{\vol(G)} \int_W \sum_{j=1}^m
(-1)^{j-1} \Omega_f^j \wedge 
\omega^1_f \wedge \cdots \wedge \widehat{\omega^j_f} \wedge \cdots
\wedge \omega_f^m \wedge r^G_f(\alpha).
\end{eqnarray*}
On a point $x \in W$ where $f = 1$, i.e. $\omega_f$ is a connection,
one gets a decomposition $TW = V_x \oplus H_x$,
then both $\Omega^j_f$ and $r^G_f9\alpha)$ are exterior forms on $H_x$,
i.e., contraction with any vector in $V_x$ is zero.
Now $\Omega^j_f \wedge r^G_f(\alpha)$ is an exterior form
of degree 
$$2 + (n -m - 1)= (n-m) + 1 = \dim H_x + 1,$$
therefore,   it must vanish.
Hence the integral above concentrates near the set where  the action
of $G$ is not locally free.
Compare with the method 
Kalkman \cite{Kal1} used to prove his formula.
The dimension counting argument will be used repeatedly below.

To prove Theorem \ref{formula1}, we take a basis $\{\xi_1, \xi_2, \xi_3\}$
for $\fg$ such that $[\xi_i, \xi_j] = a(G)\epsilon_{ijk}\xi_k$.
Here $\epsilon_{ijk}$ is nonzero only if $ijk$ is a
permutation of $123$, and when that is the case, equals the sign
of the permutation. Furthermore, $a(SU(2)) = 4\pi$, $a(SO(3)) = 2 \pi$.
For $G = SU(2)$, one can take
\begin{eqnarray*}
\xi_1 = 2\pi \left( 
	\begin{array}{cc} i & 0 \\ 0 & -i \end{array} 
	\right), & 
\xi_2 = 2\pi \left( 
	\begin{array}{cc} 0 & i \\ i & 0 \end{array} 
	\right), &
\xi_3 = 2\pi \left( 
	\begin{array}{cc} 0 & -1 \\ 1 & 0 \end{array} 
	\right).
\end{eqnarray*}
For $G = SO(3)$, one can take
\begin{eqnarray*}
\xi_1 = 2\pi \left( 
	\begin{array}{ccc} 
		0 & -1 & 0 \\ 
		1 & 0 & 0 \\
		0 & 0 & 0
	\end{array} 
	\right), & 
\xi_2 = 2\pi \left( 
	\begin{array}{ccc} 
		0  & 0 & -1 \\ 
		0  & 0 & 0 \\
		1  & 0 & 0
	\end{array} 
	\right), &
\xi_3 = 2 \pi \left( 
	\begin{array}{ccc} 
		0 & 0 & 0  \\ 
		0 & 0 & -1 \\
		0 & 1 & 0 
	\end{array} 
	\right).
\end{eqnarray*}
We also take a bi-invariant metric on $G$ such that $\{\xi_1, \xi_2, \xi_3\}$ 
defines an orthonormal basis at the identity.
Then it is clear that $\vol (T) = 1$, and 
\begin{align*}
\vol (SU(2)) & = \vol (S^3(1/2\pi)) = 2\pi^2 (1/2\pi)^3 = 1/(4\pi), \\
\vol (SO(3)) & = \vol (S^3(1/\pi))/2 = \pi^2 (1/\pi)^3 = 1/\pi.
\end{align*}
Let $\omega = \omega^j\xi_j$ be a connection on $W^0$.

\begin{lemma} For $G = SU(2)$ or $SO(3)$,
let $\alpha$ be a homogeneous $D_{\mbox{\tiny $G$}}$-closed form
of total degree $\dim (W) - 4$. Then on $\partial W$, we have
$\Omega^1 \wedge r^{\mbox{\tiny $G$}}(\alpha) = 0$.
Furthermore, we have
$$\int_{\partial W} \omega^1 \wedge \omega^2 \wedge \omega^3 
	\wedge r^{\mbox{\tiny $G$}}(\alpha) 
= -\frac{1}{a(G)} \int_{\partial W}
	\omega^1 \wedge d \omega^1 \wedge r^{\mbox{\tiny $G$}}(\alpha).
$$
\end{lemma}

\begin{proof} For each point $x \in \partial W$, 
the connection gives a decomposition $T_x{\partial W} = V_x 
\oplus H_x$. Both $\Omega^1$ and $r^{\mbox{\tiny $G$}}(\alpha)$ are 
exterior forms on $H_x$, i.e., contraction with any vector
in $V_x$ is zero. 
Then $\Omega^1 \wedge r^{\mbox{\tiny $G$}}(\alpha)$
is an exterior  form on $H_x$ of degree 
$$2 + \dim (W) - 4 = \dim (\partial W) - 1 = \dim (H_x) + 2,$$
hence it vanishes.
Now $\Omega^1 = d \omega^1 + \frac{1}{2} f^1_{bc}
\omega^b \wedge  \omega^c = d \omega^1 + a(G) \omega^2
\wedge \omega^3$, so we have
$$\omega^2 \wedge \omega^3 
= \frac{1}{a(G)} (\Omega^1 - d \omega^1).$$
Therefore,
\begin{eqnarray*}
 \int_{\partial W} \omega^1 \wedge \omega^2 \wedge
   	\omega^3 \wedge r^{\mbox{\tiny $G$}}(\alpha) 
& = &  \frac{1}{a(G)} \int_{\partial W} \omega^1
 	 \wedge (\Omega^1 - d \omega^1)  
 	 \wedge r^{\mbox{\tiny $G$}}(\alpha) \\
& = &  -\frac{1}{a(G)} \int_{\partial W} \omega^1
 	 \wedge  d \omega^1  \wedge r^{\mbox{\tiny $G$}}(\alpha).
\end{eqnarray*}
\end{proof}

\begin{proof}[Proof of Theorem \ref{formula1}]
Let $\beta = (d(f\omega^1) -  (-1 + f) u ) \wedge
p (r_f^{\mbox{\tiny $G$}}(r_f^{\mbox{\tiny $G$}}(\alpha)))
 \in \Omega_{S^1}(W)$.
It is clear that $d(f\omega^1) -  (-1 + f) u $ is $D_{\mbT}$-closed.
By Theorem 2.6, 
$$D_{\mbT}p( r_f^{\mbox{\tiny $G$}}(r_f^{\mbox{\tiny $G$}}(\alpha)))
=p(D_{\mbT} r_f^{\mbox{\tiny $G$}}(r_f^{\mbox{\tiny $G$}}(\alpha)))
= p(r_f^{\mbox{\tiny $G$}} D_{\mbT} (r_f^{\mbox{\tiny $G$}}(\alpha)))
=p( r_f^{\mbox{\tiny $G$}}(r_f^{\mbox{\tiny $G$}} D_{\mbT}(\alpha))) = 0.$$
Hence $D_{\mbT}\beta = 0$.
Furthermore, near $\partial W$, since $r^{\mbG}(\alpha)$ is basic,
$r^{\mbG}(r^{\mbG}(\alpha)) = r^{\mbG}(\alpha)$. 
So near $\partial W$, we have
$$r^{\mbT}(\beta) = r^{\mbT} (d \omega^1 \wedge
r^{\mbox{\tiny $G$}}(\alpha)) = d \omega^1 \wedge r^{\mbox{\tiny 
$G$}}(\alpha).$$
On the other hand, near each $P_k$, $f \equiv 0$ and  
 $r^{\mbox{\tiny $G$}}_f(\alpha) = \alpha$, so 
$\beta =  p(\alpha) u$. Therefore,
 by Lemma 3.1 and Kalkman's formula for $\beta$, we have 
\begin{eqnarray*}
\int_{\partial W/G} r^{\mbox{\tiny $G$}}(\alpha) 
& = & \frac{1}{\vol(G)}\int_{\partial W} \omega^1 \wedge \omega^2 \wedge
   	\omega^3 \wedge r^{\mbox{\tiny $G$}}(\alpha) \\
& = &  -\frac{1}{a(G) \vol(G)} \int_{\partial W} \omega^1
 	 \wedge  d \omega^1  \wedge r^{\mbox{\tiny $G$}}(\alpha) \\
& = &  -\frac{1}{c(G)} \int_{\partial W} \omega^1
 	 \wedge r^{\mbox{\tiny $T$}}(\beta) 
    =   -\frac{1}{c(G)} \int_{\partial W/T} r^{\mbox{\tiny $T$}}(\beta) \\
& = &  -\frac{1}{c(G)} \sum_k \int_{P_k} 
	\frac{\beta u}{\epsilon(\nu_k)} 
    =   -\frac{1}{c(G)} 
    \sum_k \int_{P_k}\frac{p(\alpha) u^2}{\epsilon(\nu_k)}.
\end{eqnarray*}
Here $c(G) = a(G) \vol(G)$, $c(SU(2)) = 4\pi \cdot 1/(4\pi) = 1$,
$c(SO(3)) =  2\pi/\pi = 2$.
\end{proof}

Now let $W$ be a compact, oriented $n$-dimensional $G$-manifold with 
$\partial W = Y \times S^2$ for 
some closed oriented manifold $Y$,
such that the action of $G$ on $\partial W$ is given by
the diagonal action of a locally-free and effective action on $Y$ and 
the coadjoint action of $SU(2)$ or $SO(3)$ on $S^2 \subset \fg^*$. 
Assume that there is a $G$-equivariant map 
$\psi: W \rightarrow S^2$,
such that $\psi|_{\partial W}$ is the projection 
$\pi_2: Y \times S^2 \rightarrow S^2$. 
Using the linear coordinates $(x_1, x_2, x_3)$ in 
the basis dual to $\{ \xi_1, \xi_2, \xi_3 \}$ on $\fg^*$,
the action of $T = \{ \exp (t\xi_1): t \in {\Bbb R} \} \subset G$ on 
$S^2 = \{ x_1^2 + x_2^2 + x_3^2 = 1 \}$
is given by
$$\exp (t\xi_1) \cdot (x_1, x_2, x_3)
=(x_1, x_2 \cos (2\pi t) - x_3 \sin (2\pi t),
	x_2 \sin (2\pi t) + x_3 \cos (2\pi t))$$
for $G =SO(3)$, and
$$\exp (t\xi_1) \cdot (x_1, x_2, x_3)
=(x_1, x_2 \cos (4\pi t) - x_3 \sin (4\pi t),
	x_2 \sin (4\pi t) + x_3 \cos (4\pi t))$$
for $G =SU(2)$.
Denote by $F$ the fixed point set of the $T$-action on $W$.
For any component $P \subset F$, 
since $\psi: W \rightarrow S^2$ is equivariant,
$\psi (P)$ is a fixed point of the $T$-action on $S^2$,
i.e., $\psi (P) = (\pm 1, 0, 0 )$.
Denote by $F_+$ the set of points fixed by $T$ which
are mapped to $(1, 0, 0)$ by $\psi$.

\begin{theorem} \label{formula2} Let $W$ be as described above. 
Assume that  $\alpha$ is an equivariantly closed 
$(n - 6)$-form  on $W$, 
such that $r^{\mbG}(\alpha|_{\partial W}) = 
\pi_1^* (\alpha_0)$ 
for some differential form $\alpha_0$ on $Y/G$,
where $\pi_1^{/G}: (Y \times S^2)/G \rightarrow Y/G$ is induced 
by the projection $\pi_1: Y \times S^2 \rightarrow Y$.
Then we have
$$\int_{Y/G} \alpha_0
	= -\frac{b(G)}{c(G)} \sum_{P_k \subset F_+}
	\int_{P_k} \frac{u^3 p(\alpha)}{\epsilon(P_k)},$$
where $c(G)$ is as in Theorem \ref{formula1},
$b(G) = 2$ for $G = SU(2)$, and $b(G) = 1$ for $G = SO(3)$.
\end{theorem}

\begin{proof}
Fix a connection $1$-form $\omega$ on the principal 
bundle $Y \rightarrow Y/G$.
The pullback of $\omega$ to $Y \times S^2$,
which we still denote by $\omega$, 
is a connection for the principal bundle 
$Y \times S^2 \rightarrow 
	(Y \times S^2)/G$.
Now on $Y$, we have 
$\Omega^1 = d \omega^1 + 4\pi \omega^2 \wedge \omega^3$.
Since $\Omega^1$ is the pullback of a $2$-form on $Y/G$,
$\Omega^1 \wedge \alpha_0$ is the pullback of a 
$(n - 4)$-form on $Y/G$. 
Hence $\Omega^1 \wedge \alpha_0 = 0$, 
since $\dim (Y/G) = n -6$.
It follows that 
\begin{equation} \label{eq:?} \begin{split} 
 \int_{Y/G} \alpha_0 
& =  \frac{1}{\vol (G)} 
	\int_Y \omega^1 \wedge \omega^2 \wedge 
	\omega^3 \wedge \alpha_0 \\
& =  \frac{1}{a(G)\vol (G)} 
	\int_Y \omega^1 \wedge (\Omega^1 - d\omega^1) 
	\wedge \alpha_0 \\
& =   -\frac{1}{c(G)} \int_Y \omega^1 \wedge
	d \omega^1 \wedge \alpha_0.
\end{split}
\end{equation}
Now we endow a symplectic structure on
$S^2 = \{ x_1^2 + x_2^2 + x_3^2 = 1 \}$ by
$$v = x_1dx_2 \wedge dx_3 + x_2dx_3 \wedge dx_1 +
	x_3dx_1 \wedge dx_2.$$
Then the  action of  $T$ on $S^2$  is a Hamiltonian action
with moment map  given by $(x_1, x_2, x_3) \mapsto 2 \pi b(G) x_1 u$. 
Hence $v_{\mbT} = v + 2\pi (1 + b(G) x_1) u$ is 
$D_{\mbT}$-closed on $S^2$. 
Now consider $\psi^*(v_{\mbT})$. 
If we denote by $V$ the vector field on $W$
generated by the action of $T$, then 
\begin{eqnarray} \label{eqn:?}
r^{\mbT}(\psi^*v_{\mbT}) 
	= \psi^*(v) - \omega^1 \wedge 
	\iota_{\mbV} \psi^*(v) + 2\pi (1 + b(G) \psi^*x_1) d\omega^1.
	\end{eqnarray}
Now on $\partial W$,
$d\omega^1 \wedge d \omega^1 \wedge 
r^{\mbG}(\alpha|_{\partial W})$
must vanish, since it is the pullback of an $(n-2)$-form 
on $Y/S^1$, which has dimension $n - 4$. 
So by $(\ref{eqn:?})$ we have
$$\omega^1 \wedge d\omega^1 \wedge 
	r^{\mbT}(\psi^*(v_{\mbT})) \wedge 
	r^{\mbG}(\alpha|_{\partial W})
= \omega^1 \wedge d\omega^1 \wedge \psi^*(v) \wedge 
	r^{\mbG}(\alpha|_{\partial W}).$$
On $\partial W$, the integration of $\psi^*(v)$  on each fiber of 
$\partial W = Y \times S^2 \rightarrow Y$ is $4 \pi$.
So from equation (\ref{eq:?}), we get
\begin{eqnarray*}
 \int_{Y/G} \alpha_0 
& = &  -\frac{1}{c(G)} \int_Y \omega^1 \wedge
	d \omega^1 \wedge \alpha_0 \\
& = &  -\frac{1}{4\pi c(G)} \int_{Y \times S^2} 
	\omega^1 \wedge d \omega^1 \wedge 
	\psi^*(v) \wedge r^{\mbG}(\alpha|_{\partial W}) \\
& = &  -\frac{1}{4\pi c(G)} \int_{Y \times S^2} 
	\omega^1 \wedge d \omega^1 \wedge 
	r^{\mbT}(\psi^*(v_{\mbT})) \wedge 
	r^{\mbG}(\alpha|_{\partial W}).
\end{eqnarray*}
Let $\beta = (d(f\omega^1) - (-1 + f) u) \wedge 
	r_f^{\mbT}(\psi^*(v_{\mbT})) \wedge 
	p (r^{\mbG}_f(r^{\mbG}_f(\alpha))) \in \Omega_{\mbT}(W)$, 
then$D_{\mbT} \beta = 0$. 
Near $\partial W$, we have 
$$r^{\mbT}(\beta) = d \omega^1 \wedge 
	r^{\mbT}(\psi^*(v_{\mbT})) 
	\wedge r^{\mbG}(\alpha|_{\partial W}). $$
On each component $P_k$ of $F_+$, 
$r^{\mbG}_f(\alpha)  = \alpha$, 
$r_f^{\mbT}(\psi^*(v_{\mbT})) = \psi^*(v_{\mbT}) = \psi^*(v) +
2\pi (1 + b(G)\psi^*x_1) = 4\pi b(G) u$.
Therefore on $F_+$, $\beta = 4\pi b(G) u^2p(\alpha)|_{P_k}$. 
Similarly, on each component of $F_-$, $\beta = 0$. 
Apply Kalkman's formula then completes the proof.
\end{proof}

\section{Application to symplectic reduction}

As a corollary to their general nonabelian 
localization formula \cite{Jef-Kir1},
Jeffrey and Kirwan proved the following:

\begin{theorem} \label{thm:JK} 
(Jeffrey-Kirwan \cite{Jef-Kir2},Corollary 3.3)
For $G = SU(2)$ or $SO(3)$, let $\mu: M \rightarrow \fg^*$ be the moment map
of a Hamiltonian $G$-action on a closed 
symplectic manifold $(M, \varpi)$
for $G = SU(2)$ or $SO(3)$. 
Suppose that $G$-action on $\mu^{-1}(0)$ is locally free and effective,
so that one can obtain the symplectic reduction  $(M_0, \varpi_0)$. 
For  any $D_{\mbG}$-closed $\eta \in \Omega_{\mbG}(M)$, 
let $\eta_0 : =  r^{\mbG}(\eta|_{\mu^{-1}(0)}) \in \Omega(M_0)$. 
Then 
$$\int_{M_o} \eta_0 e^{\varpi_0}
= - \frac{b(G)}{c(G)} Res_0 \left( u^2 \sum_{P_k \subset F_+} 
e^{ \mu_{\mbT}(P_k) u } \int_{P_k} 
\frac{p(\eta) e^{\varpi}}{\epsilon(P_k)} \right),$$
where $Res_0$ denote the coefficient of $1/u$, 
$b(G)$, $c(G)$ are the constants given in Theorem \ref{formula2}.
\end{theorem}

We will give an elementary proof of this result using the following
theorem, which is derived from Theorem \ref{formula2}.

\begin{theorem} \label{thm:local}
Assume that $\mu: M \rightarrow \fg^*$ is
as in Theorem \ref{thm:JK}.  
Suppose that $\dim(M) = 2n +6$,
then for  any $D_{\mbG}$-closed $2n$-form
$\alpha \in \Omega_{\mbG}(M)$,  we have
$$\int_{M_0} r^{\mbG}(\alpha|_{\mu^{-1}(0)})
= - \frac{b(G)}{c(G)} \sum_{P_k \subset F_+} 
 \int_{P_k} 
\frac{u^3 p(\alpha)}{\epsilon(P_k)},$$
where  $F_+$ ($F_-$)
is the subset of the fixed point set $F$ of 
$T = U(1) \subset G$ consisting of those 
components on which $\mu_{\mbT} > 0$ ($< 0$).
\end{theorem}

\begin{proof}[Proof of Theorem \ref{thm:JK} by Theorem \ref{thm:local}]
By Berline-Vergne \cite{Ber-Ver2} 
and Atiyah-Bott \cite{Ati-Bot},  one can regard the moment map
$\mu: M \rightarrow \fg^*$ as an element of $(\fg^* \otimes \Omega^0(M))^G$,
so  that  $\varpi + \mu$ is $D_{\mbG}$-closed.
Assume that $\deg (\eta) = 2s \leq 2n$. (The case of $\eta$ having
odd degree is trivial.)
Consider $\alpha:=\eta (\varpi + \mu)^{n - s}/(n-s)!$.
Then $\alpha$ is $D_{\mbG}$-closed,  and 
$r^{\mbG}(\alpha|_{\mu^{-1}(0)}) = \eta_0 \varpi_0^{n-s}/(n-s)!$.
Applying Theorem~\ref{thm:local} to this $\alpha$, we get
\begin{equation} \label{eqn:C1} \begin{split}
& \int_{M_0} \eta_0 e^{\varpi_0} 
	= \int_{M_0} \eta_0 \varpi_0^{n-s}/(n-s)! \\
= & -\frac{b(G)}{c(G)} \sum_{P_k \subset F} 
	\int_{P_k} \frac{u^3 p(\eta) (\varpi + \mu_{\mbT} u)^{n-s}}
		{(n-s)! \epsilon(P_k)}
\end{split} \end{equation}
Assume that $\dim (P_k) = 2 l_k$, then the codimension of 
$P_k$ in $M$ is $2 (n - l_k + 3)$. Now we write
\begin{align*}
p(\eta) & = \sum_a p(\eta)_a u^{s-a}, && (\deg(p(\eta)_a) = 2a) \\
(\varpi + \mu_T u)^{n-s} & = 
	\sum_b \left( \begin{array}{c} n - s \\ b \end{array} \right)
		\varpi^b \mu_{\mbT}^{n-s-b} u^{n-s-b}, && \\
\frac{1}{\epsilon(P_k)} & = 
	\frac{1}{u^{n - l_k  + 3}} \sum_c \sigma_c(P_k)/u^c, 
	&& (\deg (\sigma_c(P_k)) = 2 c) 
\end{align*}
where $p(\eta)_a$ and $\sigma_c(P_k)$ are differential forms on $P_k$. 
>From equation (\ref{eqn:C1}) we get
\begin{align*}
& \int_{M_0} \eta_0 e^{\varpi_0} \\
= & - \frac{b(G)}{c(G)} \sum_{P_k \subset F_+} \sum_{a, b, c}
	\frac{u^3}{(n-s)!u^{n - l_k  + 3}}
	\left( \begin{array}{c} n - s \\ b \end{array} \right)
	u^{(s-a) + (n-s-b) - c} \\
	& \cdot 
	\int_{P_k} p(\eta)_a \varpi^b \sigma_c(P_k) \mu_{\mbT}^{n -s-b} 
	\text{\,\,\,\,\,  (nonzero if and only if $a + b + c = l_k$)}\\
= &  - \frac{b(G)}{c(G)} \sum_{P_k \subset F_+} \sum_{a+ b+ c = l_k}
	\frac{u^{n -a -b -c +3}}{(n-s)!u^{n - l_k  + 3}}
	\left( \begin{array}{c} n - s \\ b \end{array} \right)
	\int_{P_k} p(\eta)_a \varpi^b \sigma_c(P_k) \mu_{\mbT}^{n -s-b} \\
= & - \frac{b(G)}{c(G)} \sum_{P_k \subset F_+} \sum_{a+ b+ c = l_k}
\frac{1}{(n-s)!}
	\left( \begin{array}{c} n - s \\ b \end{array} \right)
	\int_{P_k} p(\eta)_a \varpi^b \sigma_c(P_k) \mu_{\mbT}^{n -s-b}.
\end{align*}
A similar computation shows that 
$$ - \frac{b(G)}{c(G)} Res_0 \left( u^2 \sum_{P_k \subset F_+} \int_{P_k} 
\frac{p(\eta)e^{\varpi + \mu_{\mbT}(P_k) u}}
{\epsilon(P_k)} \right)$$
gives the same answer. 
 Notice now that on each $P_k$, the moment map
$\mu_{\mbT}$ is constant. This then completes the proof.
\end{proof} 

To prove Theorem \ref{thm:local}, as in Jeffrey-Kirwan 
\cite{Jef-Kir1}, one can use
the following result from symplectic geometry:

\begin{proposition} \label{propA}
(Gotay \cite{Got}, Guillemin-Sternberg
\cite{Gui-Ste}, Marle \cite{Mar}) 
Assume $0$ is a regular value of $\mu$ 
(so that $\mu^{-1}(0)  $
is a smooth manifold and  $G$ acts on $\mu^{-1}(0) $
with finite stabilizers).
Then there is a neighborhood 
${\cal O} \cong \mu^{-1}(0) \times
\{ z \in \fg^*, |z| \le h \} $
$\subseteq \mu^{-1}(0) \times \fg^*$  of $\mu^{-1} (0)$ 
on which the symplectic form is given as follows. Let
$P \stackrel{def}{=} \mu^{-1}(0) \stackrel{q}{\to} M_0$
be the orbifold  principal $G$-bundle given by
the projection map $q: \mu^{-1}(0) \rightarrow M_0 
= \mu^{-1}(0)/G$,
 and let $\omega \in \Omega^1(P) \otimes \fg$
be a connection
for it. 
Let $\varpi_0$   denote the induced symplectic form 
on $M_0$.
Then if we define a 1-form
$\tau$ on ${\cal O}\subset P \times \fg^*$ by
$\tau_{p,z}  = z(\theta)$ (for $p \in P$ and $z \in \fg^*$),
  the symplectic form on ${\cal O}$
is given by
$$ \varpi = q^* \varpi_0 + d \tau.  $$
Further,  the moment map on ${\cal O} $
is given by $\mu (p, z) = z$.
\end{proposition}

\begin{proof}[Proof of Theorem~\ref{thm:local}]  
Let $W$ be the real blow-up $\widehat{M}$ of $M$ along $\mu^{-1}(0)$, 
i.e. the result of replacing $\mu^{-1}(0)$ 
by the unit normal bundle of $\mu^{-1}(0)$ in $M$. 
Then by Proposition~\ref{propA}, $W$ is a compact manifold
with boundary $\partial W = \mu^{-1}(0) \times S^2$. 
The action of $G$ on $M$ lifts to an action on $W$,
which, on $\partial W$, is given by the diagonal action
on $\mu^{-1}(0) \times S^2$. 
Similarly, let $\widehat{\fg^*}$ be the real blowup
of $\fg^*$ at $0$. 
Then the moment map $\mu: M \rightarrow \fg^*$ 
lifts to an $G$-equivariant map 
$\hat{\mu}: W \rightarrow \widehat{\fg}^*$.
Since $\hat{\fg}^*$ can be identified with 
$S^2 \times {\Bbb R}_+$, we have a natural projection 
$\pi_1: \hat{\fg}^* \rightarrow S^2$.
Consider th composition 
$\psi = \pi_1 \circ \hat{\mu}: W \rightarrow S^2$ which on 
$\partial W$ is just the projection $\mu^{-1}(0) \times
S^2 \rightarrow S^2$.  
Then  the pair $W$ and $\psi$ satisfies the conditions in 
Theorem  \ref{formula2} and hence Theorem 4.2 follows.

\end{proof}

\begin{note}
Professor Mich\`{e}le Vergne has  suggested us to find a proof
without using any normal form theorem from symplectic geometry. 
See e.g. her note on Jeffrey-Kirwan-Witten formula \cite{Ver}.
It is actually possible in our context.
Consider in general a $G$-equivariant map $\mu: M \rightarrow \fg^*$
on a $G$-manifold $M$, 
where $G = SU(2)$ or $SO(3)$,
such that $0 \in \fg^*$ is a regular value.
Then one has a trivialization of the normal bundle of 
$\mu^{-1}(0)$ by pulling back a basis of $\fg^*$.
Applying this to the moment map of a Hamitonian $SU(2)$ or
$SO(3)$-action, we can proceed as above.
We are saved the effort of finding the normal for the symplectic form
(which is not used) in Proposition 4.1.
\end{note}


\begin{thebibliography}{99}

\bibitem{Ati-Bot} M.F. Atiyah, R. Bott, {\em The moment map and equivariant
cohomology}, {\bf Topology 23} (1984), 1-28.


\bibitem{Ber-Ver1} N. Berline, M. Vergne, 
{\em Classes caractéristiques équivariantes. 
Formule de localisation en cohomologie équivariante},
{\bf C.  R. Acad. Sci. Paris Sér. I Math. 295}
 (1982), no. 9, 539--541. 


\bibitem{Ber-Ver2} N. Berline, M. Vergne, 
{\em Fourier transforms of orbits of the coadjoint 
representation}, in 
{\bf Representation theory of reductive groups},
(Park City, Utah, 1982), 53--67, Progr. Math., 40, 
Birkhäuser Boston, Boston, MA, 1983. 

\bibitem{Ber-Get-Ver} N. Berline, E. Getzler, M. Vergne, {\bf 
Heat kernels and Dirac operator}, Springer, 1992.

\bibitem{Cao-Zho} H.-D. Cao, J. Zhou, {\em Equivariant cohomology and 
wall crossing formulas in Seiberg-Witten theory}, preprint.


\bibitem{Cao-Zho2} H.-D. Cao, J. Zhou, {\em Localization formulas on manifolds with boundaries 
and wall crossing formulas in symplectic geometry}, 
preprint, August, 1997.

\bibitem{Car1} H. Cartan, {\em Notions d'alg\`{e}bre diff\'{e}rentielle;
applications aux groupes de Lie et aux vari\'{e}t\'{e}s o\`{u} 
op\`{e}re un groupe de Lie}, in {\bf Colloque de Topologie}, 15-27. C.B.R.M. 
Bruxelles, 1950.

\bibitem{Car2} H. Cartan, {\em La transgression dans un 
groupe de lie et dans un espace fibr\'{e} principal}, in 
{\bf Colloque de Topologie}, 57-71. C.B.R.M. 
Bruxelles, 1950.

\bibitem{Duf-Kum-Ver},
M. Duflo, S. Kumar, M. Vergne,
{\em Sur la cohomologie équivariante des variétés différentiables},
{\bf Astérisque No. 215} (1993).

\bibitem{Got} M.J. Gotay, {\em On coisotropic embeddings
of presymplectic manifolds}, {\bf Proc. Amer. Math. Soc. 84}
 (1982) 111-114.

\bibitem{Gui-Ste} V. Guillemin, S. Sternberg, 
{\bf Symplectic Techniques in Physics}, 
Cambridge Univ. Press (1984).

\bibitem{Har-Law} 
R. Harvey, H.B. Lawson, Jr.,
{\em  A theory of characteristic currents associated with a singular
connection}, {\bf Astérisque No. 213} (1993).

\bibitem{Hsi}
W.-Y. Hsiang, 
{\em Cohomology theory of topological transformation groups},
 Ergebnisse der Mathematik und
ihrer Grenzgebiete, Band 85. Springer-Verlag, 
New York-Heidelberg, 1975.

\bibitem{Jef-Kir1} L.C. Jeffrey, F.C. Kirwan, 
{\em Localization for nonabelian group actions}. 
{\bf Topology 34} (1995), no. 2, 291--327.

\bibitem{Jef-Kir2} L.C. Jeffrey, F.C. Kirwan, 
{\em Intersection pairings in moduli spaces of 
holomorphic bundles on a Riemann surface},
{\bf Electronic Research Announcement of A.M.S. 1} (1995), 
Issue 2, 57-71.


\bibitem{Kal1} J. Kalkman, {\em BRST model for equivariant cohomology and 
representatives for the equivariant Thom form}, 
{\bf Commun. Math. Phys. 153} (1993), 447-463. 


\bibitem{Kal2} J. Kalkman, {\em Cohomology rings of symplectic quotients},
{\bf J. Reine Angw. Math. 458} (1995), 37-52.



\bibitem{Kam-Ton} F.W. Kamber, P. Tondeur, {\bf Foliated bundles and 
characteristic classes. Lecture Notes in
Mathematics, Vol. 493}. Springer-Verlag, Berlin-New York, 1975.



\bibitem{Law} H.B. Lawson, Jr. {\em Lecture notes of two courses on
Chern-Weil theory at Stony Brook (1991, 1992)}, taken by J. Zhou.


\bibitem{Li-Liu} T.J. Li, A. Liu, {\em General wall crossing formula},
{\bf Math. Res. Letters 2}, (1995), 797-810.

\bibitem{Mar} C.-M. Marle, 
{\em Mod\`ele d'action hamiltonienne
d'un groupe de Lie sur une vari\'et\'e symplectique}, 
in: {\bf Rendiconti del Seminario Matematico,
Universit\`a e Politechnico, Torino  43}(1985), 227-251.

\bibitem{Mat-Qui} V. Mathai, D. Quillen, {\em Superconnections, Thom classes
and equivariant differential forms}, {\bf Topology 25} (1986), 85-110.



\bibitem{Oko-Tel} C. Okonek, A. Teleman, {\em Seiberg-Witten 
invariants for manifolds with $b\sb +=1$ and the
universal wall crossing formula}, {\bf Internat. J. Math. 7 }(1996), 
no. 6, 811--832.



\bibitem{Sat} I. Satake, {\em On a generalization of the 
notion of manifold}, {\bf Proc. Nat. Acad. Sci. U.S.A. 42}
 (1956), 359--363.
 
\bibitem{Ver}
M. Vergne, 
{\em A note on the Jeffrey-Kirwan-Witten localisation formula},
{\bf  Topology 35} (1996), no. 1, 243--266. 


 
\bibitem{Wit} 
E. Witten,
{\em Two-dimensional gauge theories revisited},
{\bf J. Geom. Phys. 9} (1992), no. 4, 303--368.
 
\end{thebibliography}
\end{document}